\documentclass[12pt]{amsart}
\usepackage{amsmath} 
\usepackage{amsfonts}
\usepackage{amssymb}

\newfont{\symb}{wncyr10 scaled\magstep1}

\newcommand{\w}{\omega}

\newcommand{\Sp}{{\mathbb S}}
\newcommand{\C}{{\mathbb C}}
\newcommand{\B}{{\mathbb B}}
\newcommand{\D}{{\mathbb D}}
\newcommand{\T}{{\mathbb T}}

\newcommand{\R}{{\mathbb R}}

\newcommand{\abs}[1]{\ensuremath{\vert#1\vert}}

\normalbaselines
\newtheorem{theorem}{Theorem}[section]
\newtheorem{lemma}[theorem]{Lemma}

\newtheorem{conjecture}[theorem]{Conjecture}

\begin{document}

\title{A Weighted Estimate for the Square Function on the Unit Ball in $\C^n$}

\author[Stefanie Petermichl\and Brett D. Wick]{Stefanie Petermichl$^\dagger$\and Brett D. Wick$^\ddagger$}

\address{Stefanie Petermichl\\ Department of Mathematics\\ University of Texas at Austin\\ Austin, TX  78712\\}
\email{stefanie@math.utexas.edu}
\thanks{$\dagger$Research supported in part by a National Science Foundation Grant.}

\address{Brett D. Wick\\ Department of Mathematics\\ Vanderbilt University\\ 1326 Stevenson Center\\ Nashville, TN 37240-0001}
\email{brett.d.wick@vanderbilt.edu}
\thanks{$\ddagger$Research supported in part by a National Science Foundation RTG Grant to Vanderbilt University.}

\maketitle

\begin{abstract}
We show that the Lusin area integral or the square function on the unit ball of $\C^n$, regarded as an operator in weighted space $L^2(w)$ has a linear bound in terms of the invariant $A_2$ characteristic of the weight.  We show a dimension-free estimate for the ``area-integral'' associated to the weighted $L^2(w)$ norm of the square function. We prove the equivalence of the classical and the invariant $A_2$ classes.
\end{abstract}

\section{Introduction}

Weighted inequalities for singular integral operators appear naturally in many areas of analysis.  The theory of weights is very well understood for the ``real analysis'' case.  In a fundamental paper of Hunt, Muckenhoupt, and Wheeden, \cite{HMW}, it is shown that the so-called Calder\'on--Zygmund operators in harmonic analysis are bounded on weighted $L^p$ spaces if and only if the weight satisfies the $A_p$ condition.  Once this characterization was known, it then became of interest to determine exactly \textit{how} the norms of the operators from harmonic analysis are bounded in terms of the characteristic of the weight. One seeks the smallest power $r=r(p)$ so that $\|Tf\|_{L^p(\w)}\leq CQ_p(w)^r\|f\|_{L^p(\w)}$, with $C$ an absolute constant. The best possible power for $p=2$ usually is conjectured to be $1$. Though this is only known for very specific operators. The first successful estimate is for a dyadic analog of the square function, see \cite{HTV} and the martingale transforms in \cite{W} shortly after. As for classical Calder\'on--Zygmund operators, the best bound is only known for specific operators with certain invariance properties, such as the Hilbert and Beurling transform and the Riesz transforms. See \cite{Phil}, \cite{Priesz} and \cite{PV}. Such estimates have applications in PDE, see \cite{FKP} and \cite{PV}, and have received considerable attention.

Much of the theory of harmonic analysis can be extended to the boundary of the unit ball $\B$ since it is a domain of homogeneous type.  Many of the notions make sense on the boundary of the unit ball, but some care is needed.  In studying ``complex analysis'' questions, the choice of metric plays a distinguished role, and in this context one should use the non-isotropic metric.  This provides some additional difficulties since the ``balls'' in this metric are actually elliptical, and this provides some difference in the geometry.  So, it makes sense to ask about the boundedness properties of singular integral operators whose natural domain are functions on the boundary of the unit ball, $\Sp$, for example  the Cauchy transform.  In particular the paper \cite{LR} demonstrates that the Cauchy transform is bounded on $L^p(w)$ if and only if $w\in A_p$.  One can then inquire about other operators from harmonic analysis, such as the square function

A major motivation for this note was the paper \cite{HTV}.  In this paper, the authors determined the exact  dependence of the norm of the square function (Lusin's Area Integral, g-functions, etc.) on $L^2(\mathbb{T};w)$ in terms of the invariant characteristic of the weight.  We extend the work in \cite{HTV} to the case of the unit ball and its boundary $\Sp$. It is interesting to note that our proof has an underlying dyadic idea - similar to all other proofs of optimal bounds in weighted spaces. This is so, even though the unit sphere itself lacks a simple dyadic structure. One manages to utilize the dyadic model in one real variable to obtain a result on the unit sphere, despite the differences in geometry.

\subsection*{Acknowledgment}
The authors would like to thank Texas A \& M University and the organizers of the Workshop in Analysis and Probability.  Portions of this paper were completed while the authors were in attendance at the workshop during Summer 2006.  

The authors also thank an observant and skilled referee for many detailed observations and comments.  The presentation of the paper benefited greatly.

\subsection{Definitions} 
Let $\B$ denote the unit ball in $\C^n$, i.e. $\B=\{z \in \C^n : |z|<1\}$ and let $\Sp=\{z\in \C^n: |z|=1\}$.
We write $z=(z_1,...,z_n)$ with $z_k=x_k+iy_k$ and recall that 
$$\partial_k = \frac{\partial}{\partial z_k}= \frac12 \left(\frac{\partial }{\partial x_k}-i\frac{\partial}{\partial y_k}\right)\text{ and } \bar \partial_k = \frac{\partial}{\partial \bar z_k}= \frac12 \left(\frac{\partial }{\partial x_k}+i\frac{\partial}{\partial y_k}\right).$$
The Bergman kernel on the unit ball is given by 
$$K(z,w)=\frac{n!}{\pi^n (1-\langle z,w \rangle)^{n+1}}.$$
Define 
$$g_{ij}(z)=\partial_i  \bar\partial_j\log K(z,z)=\frac{n+1}{(1-|z|^2)^2}
[(1-|z|^2)\delta_{ij} +\bar z_i z_j]$$
The Bergman metric on $\B$ is $$\beta^2(z,\xi)=\sum_{i,j}g_{ij}(z)\xi_i\bar \xi_j.$$
So the volume element associated to the Bergman metric is 
$$dg(z)=K(z,z)dV(z)= \frac{n!}{\pi^n}\frac{dV(z)}{(1-|z|^2)^{n+1}}=\frac{d\nu (z)}{(1-|z|^2)^{n+1}}.$$Here $V$ is Lebesgue measure on the ball and $\nu$ is normalized Lebesgue measure on $\B$. 

Let $g^{ij}$ be the inverse to the matrix $g_{ij}$, so 
the Laplace--Beltrami operator on the ball is given by 
$$\widetilde \Delta = 4\sum_{i,j}g^{ij}\partial_j \bar \partial_i=4\frac{(1-|z|^2)}{n+1}\sum_{i,j}[\delta_{ij}-\bar z_i z_j]\partial_j \bar \partial_i.$$
It also has a radial form, given by
$$\widetilde \Delta f(z)=\frac{(1-r^2)}{n+1}[(1-r^2)f''(r)+\frac{2n-r^2-1}{r}f'(r)]$$ if $f(z)=f(|z|)$ for $f\in C^2(\B)$.
The invariant gradient of a $C^1(\B)$ function is the vector field given by 
$$\widetilde \nabla u = 2\sum_{i,j}g^{ij}\left(\bar\partial_i u\ \partial_j+\partial_j u\ \bar\partial_i\right).$$
The Poisson kernel for $\widetilde\Delta$ is given by $P(z,\zeta)=\frac{(1-|z|^2)^n}{|1-\langle z, \zeta \rangle|^{2n}}$ and the Green's function for $\widetilde \Delta$ is given by 
$$ G(z)=\frac{n+1}{2n}\int_{|z|}^1(1-t^2)^{n-1}t^{-2n+1}dt.$$ For $n=1$ this of course becomes $\log \frac1{|z|}$.
Let $f\in L^2(\Sp)$, define the Poisson--Szeg\"o integral $$\widetilde f(z)=\int_{\Sp}P(z,\zeta)f(\zeta)d\sigma(\zeta)$$ where $\sigma$ is normalized surface measure.
For any $\zeta \in \Sp$ we define the Kor\'anyi admissible approach region (the analogue of the non-tangential approach region or cone in the unit disc) with aperture $a>0$ as 
$$\Gamma_a(\zeta) = \{z\in\B: |1-\langle z, \zeta \rangle|<a (1-|z|^2)\}.$$ 
Note that we must assume $a>\frac{1}{2}$ otherwise $\Gamma_a(\zeta)$ is empty.  It is well known that $F$ converges to $f$ admissibly almost everywhere. 
The generalized Lusin area integral or square function on the ball with respect to $\Gamma_a(\xi)$ is 
$$S_a(f)(\zeta)=\left(\int_{\Gamma_a(\zeta)}|\widetilde \nabla \widetilde f(z)|^2dg(z)\right)^{1/2}.$$

We will be concerned with weighted $L^2$ spaces on the ball. We say that a positive $w \in L^1_{loc}(\Sp)$ function is in the class $ A_2$ if 
$$ Q_2(w)=\sup_B \left(\frac1{\sigma(B)}\int_B w(\zeta) d\sigma(\zeta) \,\cdot \frac1{\sigma(B)}\int_B w^{-1}(\zeta) d\sigma(\zeta)\right)<\infty ,$$
where the supremum runs over all non-isotropic balls $B$ on $\Sp$.
Here, we will be more concerned with the invariant $A_2$ class, denoted by $\widetilde A_2$.   A weight is in $\widetilde A_2$ if and only if 
$$\widetilde Q_2(w)=\sup_{z\in\B}\widetilde w(z)\, \widetilde{w^{-1}}(z) <\infty.$$ The quantity $\widetilde Q_2(w)$ is invariant under M\"obius transforms and is therefore more suited for questions on the ball.
We will see later that these two classes are the same, i.e. that $A_2=\widetilde A_2$ and that 
$$Q_2(w) \lesssim \widetilde Q_2(w) \lesssim Q_2(w)^2 ,$$ 
where the implied constants depend upon the dimension $n$.  Here and throughout the paper $A\lesssim B$ means $A\leq CB$ for some absolute constant $C$.

Let ${\mathcal{L}}^2(w)$ denote the space of measurable functions in the ball that are square integrable with respect to the measure $\widetilde w(z)G(z)dg(z)$.  Also let $L^2(w)$ denote the space of measurable functions that are square integrable with respect to the measure $w(\xi)d\sigma(\xi)$.  We define the operator
$$\widetilde\nabla:L^2(w)\to\mathcal{L}^2(w)$$
by sending the function $f$ to $\widetilde\nabla\widetilde f$ where $\widetilde f$ is the Poisson--Szeg\"o extension.
It is an easy calculation that
$$\|S_af\|_{L^2(w)}\le c(n)\|\widetilde \nabla \widetilde f \|_{{\mathcal{L}} ^2(w)}.$$ 

For more information on some of these concepts the reader can consult \cite{G}, \cite{K}, \cite{R}, \cite{St1}, \cite{St2} or \cite{S}.

\subsection{Main Results}

The main result in this paper is the following:

\begin{theorem}
Under the assumptions above, we have:
$$\|S_a(f)\|_{L^2(w)} \lesssim \widetilde Q_2(w)\|f\|_{L^2(w)}$$ where the implied constant may depend upon the dimension $n$. Moreover, 
$$\|\widetilde \nabla(\widetilde f)\|_{{\mathcal{L}}^2(w)}\lesssim \widetilde Q_2(w)\|f\|_{L^2(w)}$$ where we have no dependence on the dimension.  Moreover, the dependence upon $\widetilde Q_2(w)$ is sharp for both inequalities.
\end{theorem}

This theorem is proved by showing that the following string of inequalities holds:
\begin{equation}
\label{string}
\|S_a(f)\|_{L^2(w)}\stackrel{\tiny(1)}{\leq} c(n)\|\widetilde \nabla \widetilde f \|_{{\mathcal{L}} ^2(w)}\stackrel{\tiny(2)}{\leq} c\,\widetilde Q_2\|f\|_{L^2(w)}.
\end{equation}
Inequality (1) is shown by changing the order of integration, while inequality (2) is demonstrated by Bellman function techniques.

The other result demonstrated in this note is the equivalence between the weight classes $A_2$ and $\widetilde A_2$.  This is a notable result since the corresponding fact fails in $\R^n$.

\begin{theorem}
With the notation above, we have the following,
$$
Q_2(w)\lesssim \widetilde Q_2(w)\lesssim Q_2(w)^2.
$$
In particular the classes $A_2$ and $\widetilde A_2$ define the same class of weights on the unit sphere $\Sp$.
\end{theorem}

This result is demonstrated in the last section.

\section{The ``Area-Integral" Estimate}

We prove part $(2)$ of inequality (\ref{string}) using Bellman functions. 

\subsection{The Bellman function and its properties}

Let us consider the following function
$$ B(X,x,w,v)=(1+\frac1{Q})X -\frac{x^2}{Qw} -\frac{Q^2x^2}{Q^2w + (4Q^2+1)w-w^2v-4Q^2/v)}$$
on the following domain
$$ \mathcal{O}=\{(X,x,w,v)>0:x^2<Xw,1\le wv \le Q\}.$$
Note that it is convenient to think of $B=Q^{-1}B_1+B_2$ where $$B_1=X-\frac{x^2}{w}$$ and 
$$B_2=X-\frac{Q^2x^2}{Q^2w + (4Q^2+1)w-w^2v-4Q^2/v)}.$$
On $\mathcal{O}$, it enjoys the following properties:
\begin{equation}\label{Bsize} 0\le B \le 2X, \end{equation}
\begin{equation}\label{Bconcavity}-d^2B\ge C\frac1{Q^2}v(dx)^2.\end{equation}
It is of course very hard to guess such a function. It was taken from \cite{HTV} where a careful analysis was done to come up with this expression. The properties stated above are a direct calculation, we refer the reader to \cite{HTV} for detail and briefly sketch here how to get the estimates. The upper estimate on $B$ for (\ref{Bsize}) is obvious.  For the positivity of $B$, we split $B$ into $B_1$ and $B_2$ and note that $X-x^2/w \ge 0$ and hence it suffices to show that 
$$\frac{Q^2x^2}{Q^2w + (4Q^2+1)w-w^2v-4Q^2/v)}\le \frac{x^2}{w}.$$ To see the latter, note that it is equivalent to 
$w^2v^2+4Q^2 \le (4Q^2+1)vw$ which in turn is equivalent to $(vw-2Q)^2\le(2Q-1)^2vw.$ The last inequality is true for 
$1\le vw\le Q$. 
To establish (\ref{Bconcavity}), observe that $-d^2B\ge \frac2{Q} vx^2(\frac{dx}{x}-\frac{dw}{w})^2$. Also notice that $B_2(X,x,w,v)=B_1(X,x,w')$ where $w'=w+\frac{(4Q^2+1)w-w^2v-4Q^2/v}{Q^2}$. Then use the chain rule to get the estimate $-d^2B_2\gtrsim \frac1{Q^2}vx^2(\frac{dw}{w})^2$. Combining the estimates on the Hessians of $B_1$ and $B_2$ gives $-d^2B\gtrsim \frac1{Q^2} v (dx)^2$.

\subsection{A Dimension-Free Littlewood--Paley Formula}

On the unit disc one has a Littlewood--Paley formula that says 
$$\frac1{2\pi}\int_{-\pi}^{\pi}f(e^{i\theta})\bar h(e^{i \theta})d\theta = \frac1{\pi} \int_{\D}
\langle \nabla \widetilde f(z), \nabla \widetilde h(z) \rangle \log \frac1{|z|}dA(z)$$ as long as $f(0)h(0)=0$. 
In \cite{Z} a Littlewood--Paley formula for the unit ball, using the Bergman metric was derived:
$$\int_{\Sp}f(\xi)\bar h(\xi) d\sigma(\xi) = C_n \int_{\B}\langle \widetilde\nabla \widetilde f(z), \widetilde\nabla \widetilde h(z) \rangle G(z) dg(z)$$
where $f,h \in L^2(\Sp)$ with $f(0)h(0)=0$. Here, the author was mostly concerned with showing that there be a finite constant $C_n$. This result is not immediate, recall that the ball, equipped with the Bergman metric is not a compact manifold as the metric blows up on $\Sp$. 

\begin{lemma}
With $u\in C^2(\B)\cap C(\overline{\B})$ and under the assumptions above, we have the dimension-free formula
$$\int_{\B}\widetilde \Delta u(z) G(z)dg(z) =\int_{\Sp}u(\xi) d\sigma(\xi) - u(0),$$ here $\sigma$ is normalized surface measure.
\end{lemma}
The constants in this lemma can be seen by testing Green's formula for the unit ball $\B$ from \cite{Z} using the radial function $f(z)=f(|z|)=|z|^2$.

\subsection{The main inequality}

\begin{lemma}
Given ${\bf x} \in \R^k$ and smooth functions $B({\bf x})$ and $\widetilde{\bf v}(z)=(\widetilde v_1(z),...,\widetilde v_k(z))$ where the $\widetilde v_s$ are Poisson extensions, then we have the following formula for the invariant Laplacian for $b(z)=B(\widetilde{\bf v}(z))$:
$$\widetilde \Delta b(z)= 4\sum_{i,j}g^{ij}(d^2B(\widetilde{\bf v}(z))\bar \partial_i \widetilde {\bf v}(z),\overline{\partial_j\widetilde{\bf v}(z)})$$

If $-d^2B({\bf x})\ge F({\bf x}) (dx_s)^2$ then we have the estimate $-\widetilde \Delta b(z)\ge F(\widetilde{\bf v}(z))|\widetilde \nabla \widetilde v_s(z)|^2$.  Here $(dx_s)^2$ is the operator represented by the matrix with an entry $1$ in the corresponding to the second derivative in the s$^{th}$ variable and $0$ entries everywhere else. 
\end{lemma}
 
\begin{proof}
The proof is certainly a direct computation, using harmonicity of the entries of $\widetilde{\bf v}$.
\end{proof}

Applied to our situation we have 
$$-\widetilde \Delta b(z)\ge C\frac1{Q^2}\widetilde{w^{-1}} |\widetilde \nabla \widetilde f|^2.$$
Notice that with our choice of $Q=\widetilde Q(w)$ and choices for variables all entries are in our domain $\mathcal{O}$.  Using equation (\ref{Bsize}) and Green's formula applied to $b$ we have:
\begin{eqnarray*}
\lefteqn{\int_{\Sp}\frac{f^2}{w}d\sigma \ge b(0) - \int_{\Sp}bd\sigma}\\
&=& \int_{\B} -\widetilde \Delta b(z) G(z) dg(z)\\
&\ge& C\frac1{Q^2}\int_{\B}|\widetilde \nabla \widetilde f(z)|^2 \widetilde{w^{-1}}(z)G(z)dg(z).     
\end{eqnarray*}

This then proves the second inequality.  This estimate is valid for any weight $w\in A_2$.  Replacing $w$ by $w^{-1}$ and recalling that $w^{-1}\in A_2$ as well, we arrive at the desired result.

\section{The Square Function Estimate}
We now turn to (1) of inequality (\ref{string}).  Let $\mathcal{L}^2(w)$ denote the space of measurable functions on the unit ball $\B$ that are square integrable with respect to $\widetilde w(z)G(z)dg(z)$.  We have
$$
\|S_a(f)\|_{L^2(w)}\leq c(n)\|\widetilde\nabla \widetilde f\|_{\mathcal{L}^2(w)}.
$$

This follows from the following observation, $\frac{1}{\sigma(Q)}\int_{Q}w(\eta)d\sigma(\eta)\leq c(n)\widetilde w(z)$. Indeed, use the observation that for any non-isotropic ball $Q=Q(\xi,\delta)$ of radius $\delta$ and center $\xi\in\Sp$ there exists a $z_\xi\in\B$ such that $\sigma(Q)=(1-\abs{z_\xi}^2)^n$ with $\xi= \frac{z_{\xi}}{|z_{\xi}|}$.  Additionally, by the triangle inequality for the non-isotropic metric (recall the triangle inequality also holds in $\B$), for any $\eta\in Q(\xi,\delta)$ 
$$
\abs{1-\langle z_\xi,\eta\rangle}^{1/2}\leq\abs{1-\langle \xi,\eta\rangle}^{1/2}+\abs{1-\langle z_\xi,\xi\rangle}^{1/2}\leq \delta + \sigma(Q)^{\frac1{2n}}.
$$

Then for any $Q=Q(\xi,\delta)$,

\begin{eqnarray*}
\frac{1}{\sigma(Q)}\int_Q w(\eta)d\sigma(\eta) & = & \frac{1}{\sigma(Q)}\int_Q \frac{(1-\abs{z_\xi}^2)^n}{\abs{1-\langle z_\xi,\eta\rangle}^{2n}}\frac{\abs{1-\langle z_\xi,\eta\rangle}^{2n}}{(1-\abs{z_\xi}^2)^n}w(\eta)d\sigma(\eta)\\
 & = & \frac{1}{\sigma(Q)}\int_Q\mathcal{P}_{z_\xi}(\eta)\frac{\abs{1-\langle z_\xi,\eta\rangle}^{2n}}{(1-\abs{z_\xi}^2)^n}w(\eta)d\sigma(\eta)\\
 & \leq & \frac{(\delta + \sigma(Q)^{\frac1{2n}})^{4n}}{\sigma(Q)^2}\int_Q\mathcal{P}_{z_\xi}(\eta)w(\eta)d\sigma(\eta)\\
 &=& c_1(n)\widetilde w(z_\xi).
\end{eqnarray*}

Since $c_l(n)\delta^{2n}\le \sigma(Q)\le c_u(n)\delta^{2n}$ (there is only a comparison, due to the fact the ``balls'' in this metric are elliptical) then notice that $$c_1(n)\le \frac{(1+c_u^{\frac1{2n}})^{4n}}{c_l(n)^2}.$$

\begin{eqnarray*}
\lefteqn{\int_{\Sp}S^2_a(f)(\zeta)w(\zeta)d\sigma(\zeta)}\\
 & = & \int_{\Sp}\int_{\Gamma_a(\zeta)}|\widetilde \nabla f(z)|^2dg(z)w(\zeta)d\sigma(\zeta)\\
 & = & \int_{\B}|\widetilde \nabla f(z)|^2\int_{\Sp}\textbf{1}_{\Gamma_a(\zeta)}(z)w(\zeta)d\sigma(\zeta)dg(z).
\end{eqnarray*}

For fixed $z\in\B$, let
$$
E(z):=\{\zeta\in\Sp:\abs{1-\langle \frac{z}{|z|},\zeta \rangle}^{1/2}<(a^{1/2}+1)(1-\abs{z}^2)^{1/2}\}.
$$

Then one sees by the triangle inequality for the non-isotropic metric that for any $z\in\Gamma_a(\zeta) \Rightarrow \zeta \in E(z)$.   Note that $E(z)$ is a non-isotropic ball with center $\frac{z}{|z|}$ and radius $(a^{1/2}+1)(1-\abs{z}^2)^{1/2}$.  Continuing our estimate we have,

\begin{eqnarray*}
\lefteqn{\int_{\Sp}S^2_a(f)(\zeta)w(\zeta)d\sigma(\zeta)}\\
 & = & \int_{\Sp}\int_{\Gamma_a(\zeta)}|\widetilde \nabla \widetilde f(z)|^2dg(z)w(\zeta)d\sigma(\zeta)\\
 & = & \int_{\B}|\widetilde \nabla \widetilde f(z)|^2\int_{\Sp}\textbf{1}_{\Gamma_a(\zeta)}(z)w(\zeta)d\sigma(\zeta)dg(z)\\
 & \leq & \int_{\B}|\widetilde \nabla \widetilde f(z)|^2\left(\frac{1}{\sigma(E(z))}\int_{E(z)}w(\zeta)d\sigma(\zeta)\right)\sigma(E(z))dg(z)\\
 & \leq & c_1(n)  \int_{\B}|\widetilde \nabla \widetilde f(z)|^2 \widetilde w(z)\sigma(E(z))dg(z)\\
 & \leq & c_1(n)c_u(n)(a^{1/2}+1)^{2n}\int_{\B}|\widetilde \nabla \widetilde f(z)|^2 \widetilde w(z)d\sigma(\zeta)(1-\abs{z}^2)^ndg(z).
 \end{eqnarray*}

We now use the inequality $(1-\abs{z}^2)^n\leq\frac{4n^2}{n+1}G(z)$ in the last integral from above, and continue the estimate:
 \begin{eqnarray*}
\int_{\Sp}S^2_a(f)(\zeta)w(\zeta)d\sigma(\zeta) & \leq & c_1(n)c_u(n)(a^{1/2}+1)^{2n} \frac{4n^2}{n+1}\int_{\B}|\widetilde \nabla \widetilde f(z)|^2 \widetilde w(z)G(z)dg(z)\\
 & = & c(n)\int_{\B}|\widetilde \nabla \widetilde f(z)|^2 \widetilde w(z)G(z)dg(z).
 \end{eqnarray*}




This proves inequality (1) in (\ref{string}).

\section{Sharpness of the linear dependence on $\widetilde Q_2(\w)$.}

To see that our estimate is sharp, we supply a family of examples. We utilize power weights adapted to the non-isotropic metric on the sphere. Let $\w_{\alpha}(\xi)=|1-\langle \xi, \eta_0 \rangle|^{\alpha}$ for some fixed $\eta_0 \in \Sp$. Choose $f_{\alpha}(\xi)=|1-\langle \xi, \eta_0 \rangle|^{-\alpha}$. One can observe that $\w_{\alpha} \in \widetilde A_2$ if and only if $|\alpha|<n-1$. To see this, one integrates using the following formula (valid for $n\geq2$) from \cite{S}:
$$
\int_{\Sp} g(\langle \xi,\eta\rangle)d\sigma(\xi) = \frac{n-1}{\pi}\int_{\T}\int_{0}^1(1-r^2)^{n-2}g(re^{i\theta}) rdr d\theta
$$
with $g(z)=|1-z|^{\alpha}$. The largest contribution of the integrand in our case appears near $\eta_0$. Recall that $\widetilde Q_2(\w_{\alpha})=\sup_z \widetilde\w_{\alpha}(z)\widetilde{\w^{-1}_{\alpha}}(z)$. The supremum here is attained for $z$ on the ray from 0 to $\eta_0$.  
One sees that for $\alpha \to n-1$ we have $\widetilde Q_2(\w_{\alpha})\sim (n-1-\alpha)^{-1}$. Similarly one calculates that $\|f_{\alpha}\|_{L^2(\w_{\alpha})} \sim (n-1-\alpha)^{-1/2}$. It remains to see that $\|\widetilde\nabla \widetilde f_\alpha\|_{\mathcal{L}^2(w)}$, $\|S_af_{\alpha}\|_{L^2(\w_{\alpha})} \gtrsim (n-1-\alpha)^{-3/2}$. Letting $\alpha \to n-1$ then establishes sharpness. To see this, one computes the invariant gradient of the Poisson extension of the function $f_{\alpha}$ directly. One then establishes the weighted norm of the area integral similar to \cite{HTV} just more computationally involved.

\section{The Comparison of Classical and Invariant $A_2$}

In this section we establish the that the class $A_2$ and $\widetilde A_2$ are in fact the same.  This is shown by demonstrating that 
$$Q_2(w)\lesssim \widetilde Q_2(w) \lesssim Q_2^2(w)$$
where the implied constants only depend upon the dimension. 

We begin with the following fact.  First, some notation.  Let
$$
f _Q=\frac{1}{\sigma(Q)}\int_Q f(\xi)d\sigma(\xi)\quad w(Q)=\int_Qw(\xi)d\sigma(\xi).
$$

Now, observe that by Cauchy-Schwarz 
\begin{equation*}
f_Q^2 w(Q) \le Q_2(w)(f^2w)(Q).
\end{equation*}
If $P\subset Q$, and applying the above to $\textbf{1}_{Q\setminus P}$ we get
\begin{equation}
\left(1-\frac{\sigma(P)}{\sigma(Q)}\right)^2\leq Q_2(w)(w(Q)-w(P))
\end{equation}
Upon rearrangement we arrive at: 
\begin{equation}\label{betterthandoubling}
w(P)\le \left(1-\frac{(1-\sigma(P)/\sigma(Q))^2}{Q_2(w)}\right)w(Q).
\end{equation}

Given any non-isotropic ball $Q$ on the sphere, choose $z_Q\in \B$ so that the center of $Q$ is $z_Q/|z_Q|$ and its volume is $(1-|z_Q|)^n$ (hence its radius is $\sim (1-|z_Q|)^{1/2}$). There is a one-to-one correspondence between $Q$'s and $z_Q$'s. $\widetilde w(z_Q)$ will approximately correspond to $w_Q$. 

We estimate
\begin{eqnarray*}
\lefteqn{\widetilde w(z_Q)\widetilde{w^{-1}}(z_Q)\ge 
\left(\int_{Q} P_{z_Q}wd\sigma \right) \left(\int_{Q} P_{z_Q}w^{-1}d\sigma \right)   }\\
&\gtrsim& \left(\frac{(1-|z_Q|^2)^{n}}{(1-|z_Q|)^{2n}}\right)^2 
\left(\int_{Q}wd\sigma \right) \left(\int_{Q} w^{-1}d\sigma \right)\\
&\gtrsim & \frac1{(1-|z_Q|)^{2n}}\left(\int_{Q} wd\sigma \right) \left(\int_{Q} w^{-1}d\sigma \right)\\
&\gtrsim& Q_2(w).
\end{eqnarray*}
Hence $Q_2(w)\lesssim \widetilde Q_2(w)$.

We now turn to the other inequality.  We have to estimate $\widetilde w(z)\widetilde{w^{-1}}(z)\lesssim Q_2^2(w)$ with implied constant independent of $w$ and $Q$.  Using the M\"obius invariance of $\widetilde Q_2$ it suffices to assume $z=(r,0,\ldots,0)$ with  $0<r<1$ arbitrary.  Now fix $r$ and let $Q$ be the non-isotropic ball with center $(1,0,\ldots,0)$ and $\sigma(Q)=(1-r)^n$. 

To obtain the upper estimate, we exhaust the sphere by enlarging $Q$. Let us denote by $cQ$ the ball with the same center and $c^n$-fold volume (hence $c^{1/2}$-fold radius). 

\begin{equation}
\int_{\Sp}P_{z}wd\sigma 
= \int_QP_{z}wd\sigma + \sum_{k\ge 1} \int_{2^kQ\setminus 2^{k-1}Q} P_{z}w d\sigma
\end{equation}

Observe that for all $\xi \in \Sp$ 
$$P_{z}(\xi)\le\frac{(1+r)^n}{(1-r)^n}.$$
Moreover, if $\xi \in 2^kQ \setminus 2^{k-1}Q$ we have the better estimate
$$P_{z}(\xi)\lesssim 2^{-2nk}\frac{(1+r)^n}{(1-r)^n}.$$
So we get 
$$\widetilde w(z) \lesssim \frac{(1+r)^n}{(1-r)^n}\sum_{k\ge 0}2^{-2kn}w(2^kQ)$$ and similarly for $\widetilde{w^{-1}}$.
We start to estimate the product
\begin{eqnarray*}
\lefteqn{\widetilde w(z)\widetilde{w^{-1}}(z) 
\lesssim \frac{(1+r)^{2n}}{(1-r)^{2n}} \sum_{k\ge 0}2^{-2nk}w(2^kQ) \cdot \sum_{l\ge 0}2^{-2nl}w^{-1}(2^lQ)  }\\
&\le& \frac{(1+r)^{2n}}{(1-r)^{2n}} \sum_{k\ge 0} 2^{-2nk}w(2^kQ)\sum_{l=0}^k 2^{-2nl}w^{-1}(2^lQ)\\
&&+\frac{(1+r)^{2n}}{(1-r)^{2n}} \sum_{k\ge 0} 2^{-2nk}w(2^kQ)\sum_{l=k}^{\infty} 2^{-2nl}w^{-1}(2^lQ).
\end{eqnarray*}
The latter sums are equivalent (with roles of $w$ and $w^{-1}$ switched) so we estimate the second sum only.
Iterating the estimate (\ref{betterthandoubling}) we estimate for $l\ge k$
$$w(2^kQ)\le \left(1-\frac{1-4^{-n}}{Q_2(w)}\right)^{l-k}w(2^lQ),$$ and using that $\sigma (2^kQ) = 2^{nk}(1-r)^n$ we get 
\begin{eqnarray*}
\lefteqn{\frac{(1+r)^{2n}}{(1-r)^{2n}} 
\sum_{k\ge 0} 2^{-2nk}w(2^kQ)\sum_{l=k}^{\infty} 2^{-2nl}w^{-1}(2^lQ)}\\
&\lesssim& Q_2(w)\sum_{k\ge 0}2^{-2nk}\sum_{s\ge 0}\left(1-\frac{1-4^{-n}}{Q_2(w)}\right)^s\\
&\lesssim&Q_2^2(w).
\end{eqnarray*}

\section{Concluding Remarks}

The result in \cite{LR} establishes that the Cauchy transform on $L^p(w)$ is bounded if and only if the weight is in $A_p$.  It is however not known how the norm of the Cauchy transform grows with respect to the characteristic of the weight.  Thus, it is an interesting task to understand the dependence of the norm of the Cauchy transform on weighted $L^2$ spaces of the unit sphere by means of Bellman functions.  On the unit disc, using Bellman function techniques, a very simple proof of the boundedness of the Cauchy transform on weighted $L^2$ spaces is given, see \cite{NT} for continuity and \cite{PW} for the sharp result in terms of the invariant characteristic.   In the case of the disc the sharp weighted bound for the square function was a first step and a tool to obtain the sharp bound for the Cauchy transform, which was a major motivation for this note.  

As for the Cauchy transform $C$, the goal is to provide a linear and sharp estimate of the form $\|C f\|_{L^2(\w)}\lesssim \widetilde Q_2(\w)\|f\|_{L^2(\w)}$ where the implied constant is independent of $f$, $\w$ and the dimension $n$. The claim is that the invariant characteristic of the weight is the correct one, so that the dependence is both linear and no dependence on the dimension occurs. Our dimension-free estimate on the ``area integral'' illustrates that this guess was a good one. The missing ingredient for the estimate for the Cauchy transform is the following formula: $|\widetilde \nabla C f|\sim|\widetilde \nabla \widetilde f|$, which holds in one dimension, but is not true in several variables. The estimate for the square function avoids this deficiency. 

\begin{conjecture}
Let $w\in A_2$ and $C$ denote the Cauchy transform in $\C^n$.  Then, is it true that
$$
\|C(f)\|_{L^2(w)}\lesssim \widetilde Q_2(w)\|f\|_{L^2(w)},
$$
where the implied constant does not depend upon the dimension?
\end{conjecture}

\end{document}